\input amstex
\documentstyle{amsppt}
\document

\magnification 1100

\def\gen{{\frak{g}}}

\def\len{\frak{l}}

\def\sen{\frak{s}}

\def\a{{\alpha}}
\def\l{{\lambda}}
\def\b{{\beta}}
\def\d{{\delta}}
\def\n{{\eta}}
\def\t{{\tau}}
\def\g{{\gamma}}

\def\o{{\omega}}
\def\eps{{\varepsilon}}

\def\O{{\Omega}}

\def\1b{{\bold 1}}
\def\ab{{\bold a}}
\def\bb{{\bold b}}

\def\cb{{\bold c}}

\def\eb{{\bold e}}
\def\fb{{\bold f}}
\def\hb{{\bold h}}

\def\kb{{\bold k}}

\def\xb{{\bold x}}

\def\Ab{{\bold A}}
\def\Bb{{\bold B}}
\def\Cb{{\bold C}}

\def\Gb{{\bold G}}
\def\GRb{{\Gb\Rb}}

\def\Kb{{\bold K}}

\def\Rb{{\bold R}}

\def\Ub{{\bold U}}

\def\Vu{{{\underline V}}}

\def\Pu{{{\underline P}}}
\def\Iu{{{\underline I}}}
\def\iu{{{\underline i}}}

\def\Qu{{{\underline Q}}}

\def\L{{\roman L}}

\def\Wedge{{{\ts\bigwedge}}}
\def\Id{{\roman{Id}}}
\def\Sum{{{\ts\sum}}}

\def\Prod{{{\ts\prod}}}
\def\Oplus{{{\ts\bigoplus}}}

\def\ad{\text{ad}}

\def\Proj{\roman{Proj}\,}

\def\Hom{\text{Hom}\,}
\def\Isom{\text{Isom}\,}

\def\res{\roman{res}}

\def\gdm{\roman{gdm}}
\def\gch{\roman{gch}}

\def\tr{\text{tr}\,}

\def\Ker{\roman{Ker}\,}

\def\Im{\text{Im}\,}

\def\GL{{\text{GL}}}

\def\AA{{\Bbb A}}

\def\CC{{\Bbb C}}
\def\DD{{\Bbb D}}

\def\NN{{\Bbb N}}

\def\ZZ{{\Bbb Z}}

\def\Cc{{\Cal C}}
\def\Dc{{\Cal D}}

\def\Ic{{\Cal I}}
\def\Kc{{\Cal K}}

\def\Pc{{\Cal P}}
\def\Qc{{\Cal Q}}
\def\Rc{{\Cal R}}

\def\Vc{{\Cal V}}

\def\and{{\quad\roman{and}\ }}

\def\ds{\displaystyle}                
\def\ts{\textstyle}                
\def\ss{\scriptstyle}                
                
\def\qed{\hfill $\sqcap \hskip-6.5pt \sqcup$}        
\overfullrule=0pt                                    

\def\lra{{{\longrightarrow}}}
\def\simto{{\,{\buildrel\sim\over\to}\,}}

\def\u1{{\underline 1}}

\def\oh{{\overline h}}

\def\aa{{{}^a}}
\def\11{{{}^1}}
\def\la{{\langle}}
\def\ra{{\rangle}}

\newdimen\Squaresize\Squaresize=14pt
\newdimen\Thickness\Thickness=0.5pt
\def\Square#1{\hbox{\vrule width\Thickness
	      \alphaox to \Squaresize{\hrule height \Thickness\vss
	      \hbox to \Squaresize{\hss#1\hss}
	      \vss\hrule height\Thickness}
	      \unskip\vrule width \Thickness}
	      \kern-\Thickness}
\def\Vsquare#1{\alphaox{\Square{$#1$}}\kern-\Thickness}

\title Perverse sheaves and quantum Grothendieck rings\endtitle
\rightheadtext{Perverse sheaves and quivers}
\abstract We define a quantum analogue of the Grothendieck ring
of finite dimensional modules of a quantum affine algebra of simply laced type
via an analogue of Lusztig's restriction functor
on perverse sheaves on a variety related to quivers.
We get also a new geometric construction of the tensor category
of finite dimensional modules of a finite dimensional simple Lie algebra
of type $A-D-E.$
\endabstract
\author M. Varagnolo and E. Vasserot\endauthor
\thanks
\endthanks
\thanks
Both authors are partially supported by EEC grant
no. ERB FMRX-CT97-0100.\endthanks
\endtopmatter
\document

\head 1. Introduction\endhead
Finite-dimensional representations of quantum affine algebras,
say $\Ub=\Ub_q(\L\gen)$,
have been studied from various viewpoints. 
However, little is known
on the decomposition factors of tensor product of simple modules. 
From Lusztig's work the direct sum of
the Grothendieck rings of affine Hecke algebras of type $A$
can be identified with the algebra of regular functions of the 
pro-unipotent group of upper triangular unipotent $\ZZ\times\ZZ$-matrices
with finite support, in such a way that simple modules are mapped to the 
dual canonical basis of $\Ub^+_q(\sen\len_\infty)$.
It was observe recently that the induction product of simple modules
of affine Hecke algebra should be related to conjectural multiplicative
properties of the dual canonical basis, see \cite{NLT}.
The aim of our paper is to give a similar approach for tensor product
of simple modules for all simply laced types,
using the geometric realization of quantum affine algebras
in \cite{N2}, see also \cite{GV}, \cite{Va} for type $A$.
In order to do this we give a geometric construction of a flat deformation,
denoted by $\Gb\Rb$, of the Grothendieck ring of
$\Ub$ in terms of perverse sheaves on a singular variety related to quivers.
The product is defined via an analogue of Lusztig's restriction functor.
It is not commutative in general, and $\Gb\Rb$ affords a canonical basis. 
Note that $\Gb\Rb$ and its canonical basis appeared already in \cite{N3} 
in a different form. It was also observed, there, that the elements 
of the canonical basis could be identified to simple $\Ub$-modules
with a prescribed, conjectural, filtration.
However, the construction in \cite{N3} does not give
the positivity statement in Theorem 4.3.
There is no geometric construction of the tensor category of 
finite-dimensional $\Ub$-modules. 
The positivity in Theorem 4.3
suggests that a large number of information 
on tensor products of $\Ub$-modules
can be captured from the ring $\Gb\Rb$.
In particular, we formulate a generalization of a conjecture of 
Berenstein-Zelevinsky, see \cite{BZ}.

A similar construction gives a new geometric interpretation of
the tensor category of finite dimensional $\gen$-modules, see \S 5.
It would be interesting to relate it with the tensor category 
of perverse sheaves on the affine Grassmanian of the Langlands dual group.
This question appeared independentely in \cite{M}.

This work was exposed at the `Schur memorial conference' 
at the Weizmann Institute in December 2000. 
The authors would like to thank the organizers for their hospitality.

\vskip3mm

\head 2. The Grothendieck rings\endhead
\subhead 2.1\endsubhead
Let $\gen$ be a simple complex Lie algebra
with Cartan matrix $A=(a_{ij})_{i,j\in I}$.
Let $d_i\in\{1,2,3\}$ be the coprime
positive integers such that the matrix with entries 
$b_{ij}=d_ia_{ij}$ is symmetric.
Let $\a_i$ and $\o_i$ be the simple roots 
and the fundamental weights of $\gen$. 
Set $Q=\bigoplus_{i\in I}\ZZ\a_i$, $P=\bigoplus_{i\in I}\ZZ\omega_i$ and
let $P^+, Q^+$ be the semi-groups generated by $\{\a_i\}$ and $\{\o_i\}$.
Recall that $Q$ is embedded in $P$ by the linear map such that
$\a_i\mapsto\sum_ja_{ji}\o_j$.
For any $\l\in P$, $\a,\b\in Q$ we write 
$\b\geq\a$ if $\b-\a\in Q^+$ 
(resp. we write $\l\geq\a$ if $\l-\a\in P^+$).
If $\l\in P^+$ let $V(\l)$ be the simple $\gen$-module with highest weight $\l$.
If $\l\in P$ and $V$ is an integrable $\gen$-module, 
let $V_\l\subseteq V$ be the corresponding weight subspace in $V$. 
We put $$\Wedge(\l)=\bigl\{\a\in Q^+\,\big|\,V(\l)_{\l-\a}\neq\{0\}\bigr\},\quad
\Wedge^+(\l)=\bigl\{\a\in\Wedge(\l)\,\big|\,\l\geq\a\bigr\}.$$
Let $\Rb(\gen)$ be the ring of finite dimensional representations of $\gen$.
In this paper, except in \S 5.2, 
we consider only simply laced Lie algebras.

\subhead 2.2\endsubhead
The quantum loop algebra associated to ${\frak g}$ is the
$\CC(q)$-algebra $\Ub$ generated by 
$\xb^\pm_{ir},\,\kb^\pm_{i,\pm s},\,\kb_i^{\pm 1}=\kb^\pm_{i0}$ 
$(i\in I,\,r\in\ZZ,\,s\in\NN)$
modulo the following defining relations 
$$\kb_i\kb_i^{-1}=1=\kb^{-1}_i\kb_i,
\quad [\kb^\pm_{i,\pm r},\kb^\eps_{j,\eps s}]=0,$$
$$\kb_i\xb^\pm_{jr}\kb_i^{-1}=q^{\pm a_{ij}}\xb^\pm_{jr}$$
$$(z-q^{\pm a_{ji}}w)\,\kb^\eps_j(z)\,\xb^\pm_i(w)=
(q^{\pm a_{ji}}z-w)\,\xb^\pm_i(w)\,\kb^\eps_j(z)$$
$$(z-q^{\pm a_{ij}}w)\xb^\pm_i(z)\xb^\pm_j(w)=
(q^{\pm a_{ij}}z-w)\xb^\pm_j(w)\xb^\pm_i(z)$$
$$[\xb^+_{ir},\xb^-_{js}]=\delta_{ij}{{\kb_{i,r+s}^+-\kb_{i,r+s}^-}
\over q-q^{-1}}$$
$$\sum_w\sum_{p=0}^{1-a_{ij}}(-1)^{^p}
\left[\matrix 1-a_{ij}\cr p\endmatrix\right]_i
\xb^\pm_{ir_{w(1)}}\xb^\pm_{ir_{w(2)}}\cdots\xb^\pm_
{ir_{w(p)}}\xb^\pm_{js}\xb^\pm_{ir_{w(p+1)}}\cdots\xb^\pm_{ir_{w(1-a_{ij})}}=0$$
where $i\neq j,$ $r_1,...,r_{1-a_{ij}}\in\ZZ$ and $w\in S_{1-a_{ij}}$.
Here we have set $\eps=+$ or $-$,
$$[n]=(q^{n}-q^{-n})/(q-q^{-1}),\quad
[n]!=[n][n-1]...[2],\quad
\left[\matrix m\cr p\endmatrix\right]=
{[m]!\over[p]![m-p]!},$$
$$\kb^\pm_i(z)=\sum_{r\geq 0}\kb^\pm_{i,\pm r}z^{\mp r}\and
\xb_i^\pm(z)=\sum_{r\in\ZZ}\xb_{ir}^\pm\,z^{\mp r}.$$
Let $\Ub^\pm\subset\Ub$ be the subalgebra generated by the elements 
$\xb_{i,r}^\pm$ with $i\in I$, $r\in\ZZ$.
For a future use, we also introduce the elements 
$\hb_{is}\in\Ub$, $s\neq 0$, such that
$$\kb_i^\pm(z)=\kb_i^{\pm 1}
\exp\Bigl(\pm(q-q^{-1})\Sum_{s\geq 1}\hb_{i,\pm s}z^{\mp s}\Bigr).$$
Let $\Delta$ be the coproduct defined in terms of the Kac-Moody generators
$\eb_i,\fb_i,\kb_i^{\pm 1}$, $i\in I\cup\{0\}$, of $\Ub$ as follows
$$\Delta(\eb_i)=\eb_i\otimes 1+\kb_i\otimes\eb_i,\quad 
\Delta(\fb_i)=\fb_i\otimes\kb_i^{-1}+1\otimes\fb_i,\quad
\Delta(\kb_i)=\kb_i\otimes\kb_i.$$

\subhead 2.3\endsubhead
Fix $\l=\sum_i\ell_i\o_i\in P^+$, $\a=\sum_ia_i\a_i\in Q^+$,
$G_\a=\prod_i\GL_{a_i}$, $G_\l=\prod_i\GL_{\ell_i}$.
For any algebraic group $G$ let $G^\vee$ be the set of cocharacters of $G$,
and let $G^{\vee,\ad}$ be the set of conjugacy classes in $G^\vee$.
Thus $(\CC^\times)^\vee=\{q^k; k\in\ZZ\}$.
The direct sum $\GL_m^\vee\times\GL_n^\vee\to\GL_{m+n}^\vee$
gives a semigroup structure on the sets
$X^+=\bigsqcup_{\l\in P^+}G^{\vee,\ad}_\l$, 
$Y^+=\bigsqcup_{\a\in Q^+}G^{\vee,\ad}_\a.$
The Abelian groups $X, Y$ associated to $X^+,Y^+$, are identified with 
the groups $\ZZ[q^{-1},q]\otimes_\ZZ P,$ $\ZZ[q^{-1},q]\otimes_\ZZ Q$ 
via the maps 
$$\g\mapsto\sum_i(\tr\g_i)\otimes\o_i,\quad
\n\mapsto\sum_i(\tr\n_i)\otimes\a_i,$$
where $\g_i$, $\n_i$ are the $i$-th components
of the elements $\g\in G^{\vee,\ad}_\l$, $\n\in G^{\vee,\ad}_\a$.
Hereafter we may omit the symbol $\otimes$ and write
simply $q^n\l$ instead of $q^n\otimes\l$.
Consider the $\ZZ[q,q^{-1}]$-linear map
$$\O\,:\,Y\to X,\,\a_i\mapsto [2]\o_i-\Sum_{a_{ij}=-1}\o_j.$$
Hereafter, let $\g+\n$ denote the element $\g+\O(\n)\in X$.
We write $\n\succeq\d$ if $\n,\d\in Y$ are such that $\n-\d\in Y^+$
(resp. we write $\g\succeq\n$ if $\g\in X$, $\n\in Y$ 
are such that $\g-\n\in X^+$).
We have $q^{-1}\O=A+q^{-1}B+q^{-2}A$, where 
$A$, $B$ are $\ZZ[q,q^{-1}]$-linear operators such that
$A(\a_i)=\o_i$ for all $i\in I$. 
Let 
$$(\ |\ )\,:\,(\ZZ((q^{-1}))\otimes_\ZZ Q)\times X\to\ZZ((q^{-1}))$$ 
be the $\ZZ((q^{-1}))$-bilinear form such that
$(\a_i|\o_j)=\delta_{ij}$.
Let $\O^{-1}\,:\,X\to\ZZ((q^{-1}))\otimes Q$ be the inverse of $\Omega$. 
For any $\g,\g'\in X$, we put
$$\eps_{\g\g'}=\bigl(q^{-1}\O^{-1}(\bar{\g})\,\big|\,\g'\bigr)_0,\qquad
\la\g,\g'\ra=\eps_{\g\g'}-\eps_{\g'\g},$$
where $f_0$ is the constant term of a formal series $f$,
and $\,\bar{ }\,$ is the $\ZZ$-linear involution such that $\bar q=q^{-1}$.
It is easy to see that 
$$\eps_{\g+\g',\g''}+\eps_{\g\g'}=\eps_{\g,\g'+\g''}+\eps_{\g'\g''},$$
for all $\g,\g',\g''\in X^+$.
Put $\AA=\ZZ[v,v^{-1}]$. 
Let $\Ab_X$ be the $\AA$-algebra linearly spanned by elements $e^\g$, 
$\g\in X$, such that
$$e^{\g}\cdot e^{\g'}=v^{\la\g,\g'\ra}e^{\g+\g'}.\leqno(1)$$

\subhead 2.4\endsubhead
The simple finite dimensional $\Ub$-modules
are labelled by $I$-uples of monic polynomials in $\CC(q)[t]$
with nonzero constant terms, called the Drinfeld polynomials.
If $\g=\sum_k\g_k\in X^+$ with $\g_k=q^{n_k}\o_{i_k}$ and $n_k\in\ZZ$,
let $V(\g)$ be the simple finite dimensional $\Ub$-module 
whose $i$-th Drinfeld polynomial is $P_\g^{(i)}(z)=\prod_{i_k=i}(z-q^{n_k})$.
For any $\Ub$-module $V$ and any $I$-uple of formal series 
$\psi^\pm=(\psi_i^\pm)\in\CC(q)[[z^{\mp 1}]]^I$, set
$V_\psi=\bigcup_N\bigcap_i\Ker\bigl(\kb_i^\pm(z)-\psi_i^\pm\Id\bigr)^N\subseteq V.$
If $\g=\g^+-\g^-$ with $\g^\pm\in X^+$,
we put $V_\g=V_\psi$ where $\psi_i^\pm$ is the expansion at 
$\infty$ or $0$ of the rational function
$$q^{(\g(1)|\a_i)}\cdot P_\g^{(i)}(1/qz)\cdot P_\g^{(i)}(q/z)^{-1},$$
and $P_\g^{(i)}=P_{\g^+}^{(i)}/P_{\g^-}^{(i)}.$
Let $\Cc_q$ be the category of pairs $(V,F)$ where 
$V$ is a finite dimensional $\Ub$-module such that
$V=\bigoplus_\g V_\g$,
and $F$ is a decreasing $\ZZ$-filtration
on $V$ compatible with the weight decomposition,
i.e. $F_\ell V=\bigoplus_\g(V_\g\cap F_\ell V)$ for all $\ell$.
Let $\Kb(\Cc_q)$ be the $\AA$-module with one generator for each
$(V,F)\in\roman{Ob}(\Cc_q)$ modulo the relations

\itemitem{--}
$(V,F)=(V',F')+(V'',F'')$ if there is an exact sequence
$$0\to(V',F')\to (V,F)\to(V'',F'')\to 0,$$ 

\itemitem{--}
$(V,F)=v(V',F')$ if $(V,F)$ is isomorphic to $(V',F'[1])$. 

\noindent
Fix an element $\g=\sum_{k=1}^\ell\g_k\in X^+$,
with $\g_k=q^{n_k}\o_{i_k}$,
such that $n_1\geq n_2\geq\cdots\geq n_\ell$.
The $\Ub$-module 
$V(\g_1)\otimes V(\g_2)\otimes\cdots\otimes V(\g_\ell)$
does not depend on the choice of such a decomposition of $\g$, since
$$V(q^n\o_i)\otimes V(q^n\o_j)\simeq V(q^n\o_j)\otimes V(q^n\o_i)$$
for all $n,i,j$, see \cite{Ka}.
Let denote it by $W(\g)$.
Fix a highest weight vector $v_\g\in V(\g)$.
It is known that $W(\g)$ is a cyclic $\Ub^-$-module generated by the 
monomial $w_{\g}=v_{\g_1}\otimes v_{\g_2}\otimes\cdots\otimes v_{\g_\ell}$,
see \cite{Ka}, \cite{VV}.
The geometric construction in \cite{N2} implies that
$$W(\g)=\Oplus_{\g'}W(\g)_{\g'},\quad
\xb^-_{ir}\bigl(W(\g)_{\g'}\bigr)\subseteq\Oplus_{\g''}W(\g)_{\g''},$$
where the sum is over all elements $\g''\in\g'-q^\ZZ\a_i$,
see also \cite{FM}.
Note that the element $\xb_{ir}^-$ is not 
homogeneous for the weight decomposition above.
Let $x_{ir}^{(t)}$ be its component in
$$\Oplus_{\g'}\Hom\bigl(W(\g)_{\g'},W(\g)_{\g'-q^t\a_i}\bigr).$$
Set also
$$\phi_{ir}^{(t)}=\sum_{s=0}^r\left(\matrix r\cr s\endmatrix\right)
(-1)^{r-s}q^{-st}x^{(t)}_{is}.$$
We endow $W(\g)$ with the decreasing $\ZZ$-filtration such that

\itemitem{--}
$\{0\}=F_1W(\g)_\g\subset F_0W(\g)_\g=W(\g)_\g$,

\itemitem{--}
$F_kW(\g)_{\g''}=
\sum_{i,r,t}\phi_{ir}^{(t)}\bigl(F_\ell W(\g)_{\g'}\bigr),$

\noindent
where $\ell,\g'$ are such that
$$\g'=\g''+q^t\a_i,\quad\ell=k-2r-1-g''_{i,t+1}+g''_{i,t-1},
\quad\g''=\sum_{i,k}g''_{ik}q^k\cdot\o_i.$$
We have $\bigl(W(\g),F\bigr)\in\roman{Ob}(\Cc_q).$
There is a unique surjective homomorphism of $\Ub$-modules 
$W(\g)\to V(\g)$ such that $w_\g\mapsto v_\g.$
The module $V(\g)$ is endowed with the quotient filtration.
Hereafter, the classes of the pairs 
$\bigl(V(\g),F\bigr)$, $\bigl(W(\g),F\bigr)$ 
in $\Kb(\Cc_q)$ are simply denoted by $V(\g)$, $W(\g)$.
Let $\GRb\subset\Kb(\Cc_q)$ be the $\AA$-submodule
spanned by the elements $V(\g)$.
The tensor product of two objects in $\Cc_q$
is endowed with the filtration such that
$$F_k\bigl(V_\g\otimes V'_{\g'}\bigr)=\sum_{\ell+\ell'=k+\la\g,\g'\ra}
F_\ell V_\g\otimes F_{\ell'}V'_{\g'}.\leqno(1)$$
Put
$$\gdm(V_\g,F)=\sum_\ell\dim(Gr^F_\ell V_\g)\cdot v^{\ell},\quad
\gch(V,F)=\Sum_\g\gdm(V_\g,F)\cdot e^\g,$$ 
where $Gr^F$ is the associated graded space.

\proclaim{Proposition}
(a) The map $\gch\,:\,\Kb(\Cc_q)\to\Ab_X$ is a ring homomorphism.

(b) The map $\gch\,:\,\GRb\to\Ab_X$ is injective.

(c) $\GRb$ is a subring of $\Kb(\Cc_q)$.
\endproclaim

\noindent{\sl Proof:}
Put $\g=\g^+-\g^-$, where $\g^\pm=\sum_kq^{n^\pm_k}\o_{i_k}\in X^+$.
The eigenvalue of $\hb_{ir}$ on $V_\g$ is 
$r^{-1}[r]\sum_{i_k=i}\bigl(q^{rn^+_k}-q^{rn^-_k}\bigr)$, 
see \cite{FM, (2.11)} for instance.
It is known that $\Delta(\hb_{ir})=\hb_{ir}\otimes 1+1\otimes\hb_{ir}$ modulo
the linear span of elements $m_0m^-m^0\otimes n^0n^+$, where
$m_0$ (resp. $m^-, m^0,n^0,n^+$) is a monomial in the generators
$\kb_i^\pm$ (resp. $\xb_{is}^-,\hb_{is},\hb_{is},\xb^+_{is}$)
such that $m^-,n^+$ have a non-zero degree, see \cite{D}.
Thus, the weight $\g$ subspace in $V'\otimes V''$ is 
$$\bigoplus_{\g=\g'+\g''}\bigl(V'_{\g'}\otimes V''_{\g''}\bigr).$$
Then, Claim $(a)$ follows from $(2.3.1)$ and $(2.4.1)$.
Claim $(b)$ is obvious.
Claim $(c)$ is proved in Theorem 4.3. 
\qed

\vskip3mm

\noindent For a future use we introduce the following sets
$$\Wedge(\g)=\bigl\{\n\in Y^+\,\big|\,
W(\g)_{\g-\n}\neq\{0\}\bigr\},\qquad
\Wedge^+(\g)=\bigl\{\n\in\Wedge(\g)\,\big|\,\g\succeq\n\bigr\}.$$

\vskip3mm

\noindent{\bf Remark.}
The map $\gch$ appeared first in \cite{N3}.
By \cite{Vr}, the same construction holds for Yangians.
The specialization of $\gch$ at $v=1$ first appeared in \cite{Kn},
for Yangians. The case of quantum affine algebras was done in \cite{FR}.

\vskip3mm

\noindent{\bf Example.}
We give a few computations in the case $\gen=\sen\len_2$.
To simplify we omit $\o_1$ : we write $q^n$ instead of $q^n\o_1$.
We get
$$e^{q^m}\cdot e^{q^n}=v^t\, e^{q^m+q^n}\in\Ab_X,$$
where $t=0$ if $n-m$ is zero or odd, 
and $t=(-1)^\ell$ if $n-m=2\ell$ with $\ell<0$.
We have 
$$\gch V(q^n)=e^{q^n}+e^{-q^{n+2}}\in\Ab_X,$$
and
$$W(q^n+q^{n-2})=v^{-\la q^n,q^{n-2}\ra}V(q^n)\otimes V(q^{n-2}),\quad
W(kq^n)=V(q^n)^{\otimes k}\in\Gb\Rb.$$ 
Thus,
$$\gch W(q^n+q^{n-2})=e^{q^n+q^{n-2}}+e^{q^{n-2}-q^{n+2}}+
e^{-q^{n}-q^{n+2}}+v,$$
$$\gch W(kq^n)=\sum_{i=0}^k\left[\matrix k\cr i\endmatrix\right]_v
e^{iq^n-(k-i)q^{n+2}},$$
where $\left[\matrix k\cr i\endmatrix\right]_v$ is the $v$-binomial coefficient.
Note that our normalizations are different from \cite{N3} (we use $v=t^{-1}$).

\vskip3mm

\head 3. Reminder on quiver varieties\endhead
\subhead 3.1\endsubhead
Consider the graph such that : $I$ is the set of vertices,
and there are $2\delta_{ij}-a_{ij}$ egdes between $i,j\in I$.
Each edge is endowed with the two possible orientations.
The corresponding set of arrows is denoted by $H$.
If $h\in H$ let $h'$ and $h''$
the incoming and the outcoming vertex of $h$. 
Let $\oh\in H$ denote the arrow opposite to $h$.
Fix two $I$-graded finite dimensional complex vector spaces $V,$ $W$
of graded dimension $(a_i)$, $(\ell_i)$.
Let us fix once for all the following convention :
the dimension of the graded vector space $V$
is identified with the root $\alpha=\sum_{i\in I}a_i\a_i\in Q^+$ 
while the dimension of $W$ is 
identified with the weight $\lambda=\sum_i\ell_i\omega_i\in P^+$.
Set
$$E(V,W)=\bigoplus_{h\in H}M_{\ell_{h'}a_{h''}}(\CC),\quad
L(V,W)=\bigoplus_{i\in I}M_{\ell_i,a_i}(\CC),$$
$$M_{\l\a}=E(V,V)\oplus L(W,V)\oplus L(V,W).$$
For any $(B,p,q)\in M_{\l\a}$ let $B_h$ be the component of $B$
in $\Hom(V_{h''},V_{h'})$ and set
$$m_{\l\a}(B,p,q)=\sum_h\eps(h)B_hB_\oh +pq\in L(V,V),$$
where $\eps$ is a function $\eps\,:\, H\to\CC^\times$ such
that $\eps(h)+ \eps(\oh)=0.$ 
A triple $(B,p,q)\in m_{\l\a}^{-1}(0)$ is $\spadesuit$-stable
if there is no nontrivial $B$-invariant subspace of $\Ker q$.
Let $m_{\l\a}^{-1}(0)^\spadesuit$ be the subset of $\spadesuit$-stable triples. 
The group $\CC^\times\times G_\l\times G_\a$ acts on $M_{\l\a}$ by 
$$(z,g_\l,g_\a)\cdot (B,p,q)=
(z g_\a Bg_\a^{-1},zg_\a pg_\l^{-1},z g_\l qg_\a^{-1}).$$
The action of $G_\a$ on the subset $m_{\l\a}^{-1}(0)^\spadesuit$ is free.
Consider the varieties 
$$Q_{\l\a}=\Proj\bigl(\Oplus_{n\geq 0}A_n\bigr)\and 
N_{\l\a}=m^{-1}_{\l\a}(0)/\!\!/G_\a,$$
where $/\!\!/$ is the categorical quotient, and
$$A_n=\bigl\{f\in\CC[m^{-1}_{\l\a}(0)]\,\big|
\,f\bigl(g_\a\cdot(B,p,q)\bigr)=(\det g_\a)^{-n}f(B,p,q)\bigr\}.$$ 
The variety $Q_{\l\a}$ is smooth and there is a bijection 
$Q_{\l\a}\simeq m_{\l\a}^{-1}(0)^\spadesuit/G_\a.$

\subhead 3.2\endsubhead
Let $\pi_{\l\a}\,:\, Q_{\l\a}\to N_{\l\a}$ be the affinization map. 
It is a proper map. Put $d_{\l\a}=\dim Q_{\l\a}$.
It is known that $d_{\l\a}=(\a|2\l-\a)$.
If $\a\geq\b$ the extension by zero of representations of the quiver
gives a closed embedding $N_{\l\b}\hookrightarrow N_{\l\a}$.
Set $N_\l=\bigcup_\a N_{\l\a},$ $Q_\l=\bigsqcup_\a Q_{\l\a},$
$F_\l=\bigsqcup_\a F_{\l\a},$ where $F_{\l\a}=\pi_{\l\a}^{-1}(0)$.
A triple $(B,p,q)\in m_{\l\a}^{-1}(0)$ is regular if it is $\spadesuit$-stable 
and its $G_\a$-orbit is closed.
Let $m_{\l\a}^{-1}(0)^\heartsuit\subseteq m_{\l\a}^{-1}(0)^\spadesuit$ 
be the subset of regular triples. 
Let $Q_{\l\a}^\heartsuit=m_{\l\a}^{-1}(0)^\heartsuit/G_\a$ and
$N_{\l\a}^\heartsuit=m^{-1}_{\l\a}(0)^\heartsuit/\!\!/G_\a$
be the corresponding open subsets in $Q_{\l\a}$, $N_{\l\a}$.
The map $\pi_{\l\a}$ gives an isomorphism
$Q_{\l\a}^\heartsuit\simto N_{\l\a}^\heartsuit.$
It is proved in \cite{N1}, \cite{N2} that 

\itemitem{--}
$N_{\l\a}^\heartsuit\neq\emptyset\iff\a\in\Wedge^+(\l)$,
and $Q_{\l\a}\neq\emptyset\iff\a\in\Wedge(\l),$

\itemitem{--}
$N_\l=\bigsqcup_{\a}N_{\l\a}^\heartsuit$, and
$N_{\l\b}^\heartsuit\subseteq\overline{N_{\l\a}^\heartsuit}\,\iff\,
\a\geq\b.$

\subhead 3.3\endsubhead
The fixpoint set of a bijection $\phi\,:\,X\simto X$ is denoted by $X^\phi$.
The group $\CC^\times\times G_\l$ acts on $Q_{\l\a}$, $N_{\l\a}$.
For any $k\in\ZZ$ and $(\g,\n)\in G_\l^\vee\times G_\a^\vee$ we set 
$$Q_{\g\n,k}=
\bigl(G_\a\cdot m_{\l\a}^{-1}(0)^{\spadesuit,(q^k,\g,\n)}\bigr)/G_\a,
\quad Q_{\g,k}=Q_\l^{(q^k,\g)},\quad N_{\g,k}=N_\l^{(q^k,\g)}.$$
It is known that $Q_{\g\n,k}$ 
is either empty or a connected component of $Q_{\g,k}$.
Let $\pi_{\g,k}\,:\,Q_{\g,k}\to N_{\g,k}$ be the restriction of the map $\pi_\l$. 
We set $F_{\g,k}=\pi_{\g,k}^{-1}(0)$, $F_{\g\n,k}=F_{\g,k}\cap Q_{\g\n,k}$,
$Q_{\g\n,k}^\heartsuit=Q_{\g\n,k}\cap Q_{\l\a}^\heartsuit$,
$N_{\g\n,k}^\heartsuit=\pi_\g\bigl(Q_{\g\n,k}^\heartsuit\bigr).$
The restriction of $\pi_{\g,k}$ to $Q_{\g\n,k}^\heartsuit$ is an 
isomorphism onto $N_{\g\n,k}^\heartsuit$.
It is proved in \cite{N2} that 

\itemitem{--}
$Q_{\g,k}=\bigsqcup_\n Q_{\g\n,k},$ 
$N_{\g,k}=\bigsqcup_\n N_{\g\n,k}^\heartsuit,$
and $Q_{\g\n,k}$ is connected (or empty),

\itemitem{--}
the set $N_{\g\n,k}^\heartsuit$ depends only on the conjugacy classes of $\g,\n$, 

\itemitem{--}
$N_{\g\n,1}^\heartsuit\neq\emptyset\iff\n\in\Wedge^+(\g)$
and $Q_{\g\n,1}\neq\emptyset\iff\n\in\Wedge(\g)$.

\noindent
To simplify, hereafter we set 
$Q_{\g\n}=Q_{\g\n,1},$ $Q_\g=Q_{\g,1},$ $N_\g=N_{\g,1}$, etc.
Put $$d_{\g\n}=\bigl(\bar\n\,\big|\,[2]\g-q\n\bigr)_0.$$
If $Q_{\g\n}\neq\emptyset$ then 
$d_{\g\n}=\dim Q_{\g\n}$, see \cite{N2, (4.1.6)}. 

\subhead 3.4\endsubhead
If an algebraic group $G$ acts on a variety $X$, and if $\phi\in G^\vee$, we put
$$X^{+\phi}=\{x\in X\,|\,\lim_{z\to 0}\phi(z)\cdot x\in X^\phi\},\quad
X^{-\phi}=\{x\in X\,|\,\lim_{z\to\infty}\phi(z)\cdot x\in X^\phi\}.$$
For any $k\in\ZZ$, $\g\in G_\l^\vee$ , $\t\in(\CC^\times\times G_\l)^\vee$
we have the commutative diagram
$$\matrix
Q_{\g,k}&{\buildrel\tilde\iota_\pm\over\hookleftarrow}&
Q_{\g,k}^{\pm\tau}
&{\buildrel\tilde\kappa_\pm\over\twoheadrightarrow}&Q_{\g,k}^\tau\cr
\downarrow&&\downarrow&&\downarrow\cr
N_{\g,k}&{\buildrel\iota_\pm\over\hookleftarrow}&
N_{\g,k}^{\pm\tau}
&{\buildrel\kappa_\pm\over\twoheadrightarrow}&N_{\g,k}^\tau,\cr
\endmatrix$$
where $\tilde\iota_\pm,\iota_\pm$ are the embeddings,
and $\tilde\kappa_\pm,\kappa_\pm$ are the obvious projections.
Since the map $\pi_{\g,k}$ is proper 
the left square is Cartesian.

\vskip3mm

\noindent{\bf Remark.}
The maps $\tilde\iota_\pm,\iota_\pm$ are closed embeddings.
We have $Q_{\g,k}^{\pm\tau}=\pi_{\g,k}^{-1}(N_{\g,k}^{\pm\tau})$ 
since $\pi_{\g,k}$ is a proper map.
Thus, it is sufficient to consider the case of $\iota_\pm$. 
From \cite{L2}, 
we can fix a finite set of generators of the ring 
$\CC[m_{\l\a}^{-1}(0)]^{G_\a}$
consisting of eigenvectors of the group 
$\tau(\CC^\times)\times\g(\CC^\times)\subset\CC^\times\times G_\l$.
These generators give a $\tau(\CC^\times)$-equivariant closed embedding of 
the variety $N_{\g,k}$ in a finite dimensional representation 
of $\tau(\CC^\times)$. 
But $X^{\pm\phi}$ is a closed subset of $X$ in the particular case
where $X$ is a representation of the one-parameter subgroup $\phi$.
Thus $N_{\g,k}^{\pm\tau}$ is a closed subset of $N_{\g,k}.$

\subhead 3.5\endsubhead
Fix $\l',\l''\in P^+$, fix $I$-graded vector spaces 
$W'$, $W''$ of dimension $\l'$, $\l''$, and
fix $\g'\in G_{\l'}^\vee$, $\g''\in G_{\l''}^\vee$. 
Put $\g=\g'+\g''$, $\l=\l'+\l''$, $W=W'\oplus W''$
and $\tau=q\cdot\Id_{W'}\oplus\Id_{W''}.$

\proclaim{Lemma 1}
(a) The direct sum of representations of the quiver gives an isomorphism
$Q_{\g',k}\times Q_{\g'',k}\simeq Q_{\g,k}^\t$, and a map 
$\phi\,:\,N_{\g',k}\times N_{\g'',k}\to N_{\g,k}^\t$.

(b) The map $\phi$ is finite, bijective and is
compatible with the stratifications.
\endproclaim

\noindent{\sl Proof:}
The first claim is well-known, see \cite{VV, Lemma 4.4} for instance.
Let $\phi\,:\,m_{\l'\a'}^{-1}(0)\times m_{\l''\a''}^{-1}(0)\to m_{\l\a}^{-1}(0)$
be the direct sum of representations of the quiver in \S 3.1.
The induced map $N_{\l'\a'}\times N_{\l''\a''}\to N_{\l\a}$ is a morphism
of algebraic varieties.
We have 
$$\phi\bigl(m_{\l'\a'}^{-1}(0)^\heartsuit\times 
m_{\l''\a''}^{-1}(0)^\heartsuit\bigr)\subset m_{\l\a}^{-1}(0)^\heartsuit,$$
since a triple $(B,p,q)\in m_{\l\a}^{-1}(0)$ is regular if and
only if it is stable and costable (i.e. there is  no proper
$B$-invariant subspace of $V$ containing $\Im p$), see \cite{L2}.
Fix $\n'\in G_{\a'}^\vee$, $\n''\in G_{\a''}^\vee$ such that $\n=\n'+\n''$.
By the first part, $\phi$ gives an isomorphism 
$N_{\g'\n',k}^\heartsuit\times N_{\g''\n'',k}^\heartsuit
\simto\bigl(N_{\g\n,k}^\heartsuit\bigr)^{\t}$.
In particular it induces a bijection 
$$N_{\g',k}\times N_{\g'',k}=
\bigsqcup_{\n',\n''}N_{\g'\n',k}^\heartsuit\times N_{\g''\n'',k}^\heartsuit
\simto 
\bigsqcup_\n\bigl(N_{\g\n,k}^\heartsuit\bigr)^\tau=N_{\g,k}^\t,$$
which is compatible with the stratifications.
This map is clearly affine, since $N_{\l\a}$ is an affine variety.
Thus it is finite.
\qed

\vskip3mm

\noindent
If $k=1$ we get
$$\matrix
Q_\g&{\buildrel\tilde\iota_\pm\over\hookleftarrow}&Q_\g^{\pm\t}
&{\buildrel\tilde\kappa_\pm\over\twoheadrightarrow}&
Q_\g^\t&\simeq&
Q_{\g'}\times Q_{\g''}\cr
{\ss\pi_\g}\downarrow&$\qed$&\downarrow&&\downarrow&&\downarrow\cr
N_\g&{\buildrel\iota_\pm\over\hookleftarrow}&
N_\g^{\pm\t}&{\buildrel\kappa_\pm\over\twoheadrightarrow}&N_\g^\t
&{\buildrel\phi\over\leftarrow}&N_{\g'}\times N_{\g''}.
\endmatrix$$
Fix $\n'\in G^\vee_{\a'},\n''\in G^\vee_{\a''}$.
Let $\kappa^\pm_{\n'\n''}$ be the relative dimension of the map $\tilde\kappa_\pm$
above the component $Q_{\g'\n'}\times Q_{\g''\n''}$.
Set $\n=\n'+\n''$.

\proclaim{Lemma 2} We have

(a) $\kappa^+_{\n'\n''}+\kappa^-_{\n'\n''}=d_{\g\n}-d_{\g'\n'}-d_{\g''\n''},$

(b) $\kappa^\pm_{\n'\n''}=\kappa^\mp_{\n''\n'},$

(c) If $\d'\in\Wedge^+(\g')$, $\d''\in\Wedge^+(\g'')$ are such that 
$\n'\succeq\d'$, $\n''\succeq\d''$, then 
$$\eps_{\g'\g''}-\eps_{\g'-\d',\g''-\d''}=\kappa^\pm_{\n'\n''}-
\kappa^\pm_{\n'-\d',\n''-\d''}=\kappa^\pm_{\d'\d''}.$$
\endproclaim

\noindent{\sl Proof:}
Part $(a)$ is immediate. Let us check Part $(b)$.
The one-parameter subgroup $q\cdot\Id_{W'}\oplus\Id_{W''}$
acts fiberwise on the normal bundle to 
$Q_{\g'\n'}\times Q_{\g''\n''}$ in $Q_{\g\n}$.
By definition $\kappa^\pm_{\n'\n''}$ is the dimension of the attracting
(resp. repulsing) subbundle.
The class in equivariant $K$-theory of the tangent bundle to $Q_{\g\n}$ 
is given in \cite{N1, \S 4.1}. We get 
$$\kappa^+_{\n'\n''}=
\bigl(\bar\n'\,\big|\,q^{-1}\g''\bigr)_0+
\bigl(\bar\n''\,\big|\,q\g'\bigr)_0-
\bigl(\bar\n''\,\big|\,q\O(\n')\bigr)_0,\leqno(1)$$
and $\kappa^-_{\n'\n''}=\kappa^+_{\n''\n'}.$
Observe that
$$\bigl(\O^{-1}(\g)\,\big|\,\O(\n)\bigr)=
(\n\,\big|\,\g),\qquad\forall\g\in X,\,\n\in Y.$$
Part $(c)$ is proved by a direct computation using $(3.5.1)$ and
$$\matrix
\eps_{\g'\g''}-\eps_{\g'-\d',\g''-\n''}
&=\bigl(q^{-1}\O^{-1}(\bar\g')\,\big|\,\g''\bigr)_0
-\bigl(q^{-1}\O^{-1}(\bar\g')-q^{-1}\bar\d'\,\big|\,
\g''-\O(\d'')\bigr)_0\hfill\cr
&=\bigl(q^{-1}\O^{-1}(\bar\g')\,\big|\,\O(\d'')\bigr)_0
+\bigl(q^{-1}\bar\d'\,\big|\,\g''\bigr)_0
-\bigl(q^{-1}\bar\d'\,\big|\,\O(\d'')\bigr)_0\hfill\cr
&=\bigl(q\bar\d''\,\big|\,\g'\bigr)_0
+\bigl(q^{-1}\bar\d'\,\big|\,\g''\bigr)_0
-\bigl(q\bar\d''\,\big|\,\O(\d')\bigr)_0.\hfill
\endmatrix$$
\qed

\vskip3mm

\head 4. The product\endhead
\subhead 4.1\endsubhead
For any complex algebraic variety $X$, let $\Dc(X)$ be the bounded derived 
category of complexes of constructible sheaves of $\CC$-vector spaces on $X$.
For any irreductible local system $\phi$ on a locally closed set 
$Y\subset X$, let $IC(Y,\phi)$ be the corresponding intersection 
cohomology complex. Let $\CC_Y$ be the constant sheaf on $Y$.
We set $IC(Y)=IC(Y,\CC_Y)$.
Recall that the direct image of a simple perverse sheaf 
by a finite bijective map is still a simple perverse sheaf.
Let $\DD$ denote the Verdier duality.

Fix $\g,\g'\in X^+$, $\n,\n\in Y^+$.
Fix $\l,\a$ such that $\g\in G_\l^{\vee,\ad}$, $\n\in G_\a^{\vee,\ad}$.
Hereafter we may identify a cocharacter in $G^\vee_\l$, $G^\vee_\a$,
and its conjugacy class in $X^+, Y^+$. 
Let $\Dc(N_\g)^\heartsuit$,
$\Dc(N_\g\times N_{\g'})^\heartsuit$ 
be the full subcategories of 
$\Dc(N_\g)$,
$\Dc(N_\g\times N_{\g'})$ 
consisting of all complexes which are constructible 
with respect to the stratification in $\S 3.3$. 
Set $IC_{\g\n}=IC(N_{\g\n}^\heartsuit)$,
$\CC_{\g\n}=\CC_{N_{\g\n}^\heartsuit}[d_{\g\n}]$
for any $\n\in\Wedge^+(\g)$, and
$L_{\g\n}=\pi_{\g!}\CC_{Q_{\g\n}}[d_{\g\n}]$ 
for any $\n\in\Wedge(\g)$.
Let $\Qc_\g$, $\Qc_{\g\g'}$ 
be the full subcategories of 
$\Dc(N_\g)^\heartsuit$,
$\Dc(N_\g\times N_{\g'})^\heartsuit$ 
consisting of all complexes which are isomorphic to 
finite direct sums of the sheaves $IC_{\g\n}[k]$, 
$IC_{\g\n}[k]\boxtimes IC_{\g'\n'}[k']$, $k,k'\in\ZZ$.
The complex $L_{\g\n}$ belongs to $\roman{Ob}(\Qc_\g)$, 
see \cite{N2, Theorem 14.3.2}.
If $\g',\g''$, $\iota_\pm,\kappa_\pm$, $\t$ are as in $\S 3.5$, 
we have the functor 
$$\res^\pm_{\g'\g''}=\kappa_{\pm!}\iota_\pm^*
\,:\,\Dc(N_\g)^\heartsuit\to\Dc(N_\g^\t)^\heartsuit.$$

\proclaim{Lemma}
We have

$(a)$ $\res^\pm_{\g'\g''}(L_{\g\n})=
\Oplus_{\n=\n'+\n''}\phi_!(L_{\g'\n'}\boxtimes L_{\g''\n''})
[\kappa^\mp_{\n'\n''}-\kappa^\pm_{\n'\n''}],$

$(b)$ $\DD\circ\res^\pm_{\g'\g''}=\res^\mp_{\g'\g''}\circ\DD,$
and $\res^\pm_{\g'\g''}=\res^\mp_{\g''\g'}.$

$(c)$ For any complex $P\in\roman{Ob}(\Qc_\g)$ there is a
complex $P'\in\roman{Ob}(\Qc_{\g'\g''})$ such that
$\res^\pm_{\g'\g''}(P)\simeq\phi_!(P)$.
The complex $P'$ is unique up to isomorphism.

\endproclaim

\noindent{\sl Proof:}
By base change, the diagram in \S 3.5 gives
$$\res^\pm_{\g'\g''}(L_{\g\n})=\pi_{\g !}\tilde\kappa_{\pm!}\tilde\iota_\pm^*
\CC_{Q_{\g\n}}[d_{\g\n}].$$
From \cite{L1, 8.1.6} the complex
$\pi_{\g !}\tilde\kappa_{\pm!}\tilde\iota_\pm^*\CC_{Q_{\g\n}}$
is semi-simple, and there are short exact sequences of perverse sheaves
$$0\to{}^pH^n(f_j)_!\tilde\iota_\pm^*\CC_{Q_{\g\n}}\to 
{}^pH^n(f_{\leq j})_!\tilde\iota_\pm^*\CC_{Q_{\g\n}}\to 
{}^pH^n(f_{\leq j-1})_!\tilde\iota_\pm^*\CC_{Q_{\g\n}}\to 0,$$
where ${}^pH^n$ is the perverse cohomology, and $f_j$ (resp. $f_{\leq j}$)
is the restriction of the map $\pi_\g\tilde\kappa_\pm$ to the union of 
all subvarieties 
$$\tilde\kappa_\pm^{-1}
\bigl(Q_{\g'\n'}\times Q_{\g''\n''}\bigr)\subset Q_\g^{\pm\t}$$
of dimension $j$ (resp. $\leq j$).
We have also
$$\pi_{\g!}\tilde\kappa_{\pm!}\tilde\iota_\pm^*\CC_{Q_{\g\n}}[d_{\g\n}]
=\phi_!\bigl(L_{\g'\n'}\boxtimes L_{\g''\n''}\bigr)
[d_{\g\n}-2\kappa^\pm_{\n'\n''}].$$
Thus, Claim $(a)$ follows from Lemma 3.5.2.$(a)$. 
Claim $(b)$ is due to the auto-duality of $L_{\g\n}$,
since the map $\pi_\g$ is proper, and Lemma 3.5.2.$(b)$. 
The first part of Claim $(c)$ follows from Claim $(a)$, 
since a direct summand of a complex in $\Qc_\g$ belongs to $\Qc_\g$.
The second part of Claim $(c)$ is due to Lemma 3.5.1.$(b)$.
\qed 

\vskip3mm

\subhead 4.2\endsubhead
Let $\Kc_\g$ be the $\AA$-module with one generator 
for each isomorphism class of object of $\Qc_\g$,
with relations
$P+P'=P''$ if the complex $P''$ is isomorphic to $P\oplus P'$, and
$P=vP'$ if the complex $P$ is isomorphic to $P'[1]$. 
The elements $IC_{\g\n}$, with $\n\in\Wedge^+(\g)$,
form a $\AA$-basis of $\Kc_\g$. 
Let $\res_{\g'\g''}$ be the $\AA$-linear map
$\Kc_\g\to\Kc_{\g'}\otimes\Kc_{\g''}$
such that
$$\res_{\g'\g''}(P)=v^{\la\g',\g''\ra}\sum_iP_i'\otimes P_i''$$
where
$\res^+_{\g'\g''}(P)=\Oplus_i\phi_!\bigl(P_i'\boxtimes P_i''\bigr).$
It is well-defined and unique by Lemma 4.1.$(c)$.

\proclaim{Lemma 1}
(a) In $\Kc_\g$ we have
$$L_{\g\n}=\Sum_\d\gdm V(\g-\d)_{\g-\n}\,IC_{\g\d}.$$
In particular, the elements $L_{\g\n}$, with $\n\in\Wedge^+(\g)$,
form a $\AA$-basis of $\Kc_\g$. 

(b) If $\d\in\Wedge^+(\g)$ there is a unique surjective map
$\Kc_\g\to\Kc_{\g-\d}$ such that 
$L_{\g\n}\mapsto L_{\g-\d,\n-\d}$ if 
$\n\in\Wedge(\g)$, $\n\succeq\d$, 
and $L_{\g\n}\mapsto 0$ else.

(c) If $\d\in\Wedge^+(\g)$, $\d'\in\Wedge^+(\g')$, 
$\d''\in\Wedge^+(\g'')$, $\d=\d'+\d''$, the square
$$\matrix
\Kc_\g&{\buildrel\res\over\to}&\Kc_{\g'}\otimes\Kc_{\g''}\cr
\downarrow&&\downarrow\cr
\Kc_{\g-\d}&{\buildrel\res\over\to}&\Kc_{\g'-\d'}\otimes\Kc_{\g''-\d''}
\endmatrix$$
is commutative.
\endproclaim

\noindent{\sl Proof:}
Fix $\d\in\Wedge^+(\g)\cap G_\b^{\vee,\ad}$ such that $\n\succeq\d$, 
and fix $x_\d\in N_{\g\d}^\heartsuit$. 
We have an isomorphism 
$$W(\g-\d)_{\g-\n}\simeq
\Oplus_kH_k(F_{\g-\d,\n-\d})\simeq
\Oplus_kH_k\bigl(Q_{\g\n}\cap\pi^{-1}_\g(x_\d)\bigr)$$
such that $w_{\g-\d}\in H_0(F_{\g-\d,0}),$
see \cite{VV, Theorem 7.12}, \cite{N2, Theorems 3.3.2 and 7.4.1}.
We first check that
$$\gdm W(\g-\d)_{\g-\n}=
\sum_kv^{d_{\g-\d,\n-\d}-k}\dim H_k\bigl(Q_{\g\n}\cap\pi^{-1}_{\g}(x_\d)\bigr).
\leqno(1)$$
To simplify the notations, we may assume that $\d=0$, 
without loss of generalities.
Let $C_{\l,\a+\a_i,\a}\subseteq Q_{\l,\a+\a_i}\times Q_{\l\a}$
be the set of pairs $(x',x)$ such that 
$x$ is a subrepresentation of $x'$.
For any $\n,\n'$ put 
$$C_{\n'\n}=C_{\l,\a+\a_i,\a}\cap\bigl(Q_{\g\n'}\times Q_{\g\n}\bigr).$$
If $C_{\n'\n}\neq\emptyset$ then
$\n'=\n+q^t\a_i$ for some $t\in\ZZ$.
Set 
$$d_{\n'\n}=\dim C_{\n'\n},\quad
e_{\n'\n}=d_{\g\n}+d_{\g\n'}-2d_{\n'\n}.$$
Let $\star$ be the convolution product
in Borel-Moore homology, see \cite{CG}.
By definition, we have
$$H_{d_{\g\n'}+d_{\g\n}-e}^{BM}\bigl(C_{\n'\n}\bigr)\star 
H_{d_{\g\n}-\ell}(F_{\g\n})
\subseteq H_{d_{\g\n'}-k}(F_{\g\n'}),\quad k=\ell+e,$$
see \cite{CG, Lemma 8.9.5}.
Recall that $x^{(t)}_{ir}$ acts on $H_*(F_\g)$
by the $\star$-product by an element of the form
$$\sum_\n(\theta_{\n'\n}\cup q^{rt}e^{r\o_{\n'\n}})\cap[C_{\n'\n}]\in 
H^{BM}_*(C_{\n'\n}),$$
where $\n'=\n+q^t\a_i$, $[C_{\n'\n}]$ is the fundamental class,
$\o_{\n'\n},\theta_{\n'\n}\in H^{2*}(Q_{\g\n'}\times Q_{\g\n})$,
$\deg\o_{\n'\n}=2$, and $\theta_{\n'\n}$ is invertible.
Moreover, $\o_{\n'\n}$, $\theta_{\n'\n}$ do not depend on $r$. 
More precisely, from \cite{N2, (9.3.2), \S 13.4}, we have 
$$\theta_{\n'\n}=e^{k\o_{\n'\n}}\cup(1\otimes\nu_{\n'\n})\leqno(2)$$
where $k\in\ZZ$ and $\nu_{\n'\n}\in H^{2*}(Q_{\g\n})$ is invertible.
Fix a non-zero $v\in H_0(F_{\g0})$. 
The space $H_*(F_\g)$ is spanned by the elements
$\xb_{i_1r_1}^-\cdots\xb^-_{i_sr_s}(v)$. 
Thus, for any $\n'\in Y^+\setminus 0$, we get
$$H_*(F_{\g\n'})=\sum_{i,t,r}x^{(t)}_{ir}\star H_*(F_{\g\n}),$$
where $\n=\n'-q^t\a_i$. 
Set $\psi_{ir}^{(t)}=\sum_\n\o_{\n'\n}^r\cap[C_{\n'\n}].$
Using (4.2.2) we get
$$H_*(F_{\g\n'})=\sum_{i,t,r}\psi^{(t)}_{ir}\star H_*(F_{\g\n}).$$
The $\star$-product by $\psi_{ir}^{(t)}$ on $H_*(F_{\g\n})$
is a homogeneous operator of degree
$e_{\n'\n}+2r\in\ZZ$. Thus,
$$H_{d_{\g\n}-k}(F_{\g\n'})=\sum_{i,t,r}\psi^{(t)}_{ir}
\star H_{d_{\g\n}-\ell}(F_{\g\n}),$$
where $k=e_{\n'\n}+\ell+2r$.
Set
$$F_\ell H_*(F_{\g\n})=\Oplus_{\ell'\geq\ell}H_{d_{\g\n}-\ell'}(F_{\g\n}).$$
A direct computation gives
$$\phi_{ir}^{(t)}=\sum_\n\theta_{\n'\n}\cap(e^{\o_{\n'\n}}-1)^r\cap[C_{\n'\n}],$$
where $\n'=\n+q^t\a_i$. 
Thus
$$F_kH_*(F_{\g\n'})=\sum_{i,t,r}\phi^{(t)}_{ir}
\star F_\ell H_*(F_{\g\n}),\leqno(3)$$
where $k=e_{\n'\n}+\ell+2r$.

The $\g$-fixed part of the complex 
\cite{N2, (5.1.1)} is the normal bundle
of $C_{\n'\n}$ in $Q_{\g\n'}\times Q_{\g\n}$. Thus
$$d_{\g\n}+d_{\g\n'}-d_{\n'\n}=\bigl(q\bar\n+q^{-1}\bar\n'\,\big|\,\g\bigr)_0
-\bigl(q\bar\n\,\big|\,\O(\n')\bigr)_0.$$
From $\n'=\n+q^t\a_i$, we get 
$$e_{\n'\n}=\bigl(\bar\n'-\bar\n\,\big|\,q^{-1}(\g-\n)-q(\g-\n')\bigr)_0
=1+\bigl(\a_i\,\big|\,q^{-t}(q^{-1}-q)(\g-\n')\bigr)_0.$$
Using (4.2.3) and \S 2.4 we get
$$F_\ell W(\g)_{\g-\n}=F_\ell H_*(F_{\g\n}).$$
The identity $(4.2.1)$ follows.

To prove Lemma 4.2.1.$(a)$ set
$L_{\g\n}=\Oplus_{k,\d\preceq\n}M_{\d k}\otimes IC_{\g\d}[k]$.
If $\g,\n\succeq\d$, let
$$\phi_{\d k}\,:\,H_{d_{\g-\d,\n-\d}-k}(F_{\g-\d,\n-\d})\to
H^{d_{\g-\d,\n-\d}+k}(F_{\g-\d,\n-\d}),$$
be the composition of the chain of maps 
$$H_{*-k}(F_{\g-\d,\n-\d})\to
H^{BM}_{*-k}(Q_{\g-\d,\n-\d})\to
H^{*+k}(Q_{\g-\d,\n-\d})\to
H^{*+k}(F_{\g-\d,\n-\d}).$$
A detailed analysis of the gradings in \cite{N2, \S 14}, \cite{CG, \S 8}
shows that $M_{\d k}=\Im\phi_{\d k}.$ 
Since $V(\g-\d)_{\g-\n}\simeq\Oplus_k M_{\d k},$
we get
$\gdm V(\g-\d)_{\g-\n}=\sum_kv^k\dim M_{\d k}.$

Let us prove part $(b)$.
By \cite{N2, Theorem 3.3.2} we have for any $\delta\in\Wedge^+(\g)$
$$N^\heartsuit_{\g-\d,\n-\d}=\emptyset\iff
N^\heartsuit_{\g\n}=\emptyset,\quad
Q_{\g-\d,\n-\d}=\emptyset\iff
Q_{\g\n}=\emptyset.$$
Thus, using \S 3.3 we get
$$\matrix
\n-\d\in\Wedge^+(\g-\d)&\iff&
\n\in\Wedge^+(\g),\quad\n\succeq\d,\cr
\n-\d\in\Wedge(\g-\d)&\iff&
\n\in\Wedge(\g),\quad\n\succeq\d.
\endmatrix\leqno(4)$$
By (4.2.4) there is a unique surjective map
$\Kc_\g\to\Kc_{\g-\d}$ such that 
$$IC_{\g\n}\mapsto IC_{\g\n}\quad\roman{if}\quad\n\succeq\d,
\quad\roman{and}\quad IC_{\g\n}\mapsto 0\quad\roman{else.}$$
Using (4.2.4) again and Claim $(a)$ of the lemma, we see that this map
satisfies the requirements in Claim $(b)$.

Set 
$$\matrix
A=\kappa^-_{\n'\n''}-\kappa^+_{\n'\n''}+\la\g',\g''\ra,\hfill\cr
B=\kappa^-_{\n'-\d',\n''-\d''}-\kappa^+_{\n'-\d',\n''-\d''}+
\la\g'-\d',\g''-\d''\ra.\hfill
\endmatrix$$
Using Lemma 3.5.2.$(b),(c)$ we get $A=B$.
Thus, Claim $(c)$ follows from Claim $(b)$ and Lemma 4.1.$(a)$.
\qed

\vskip3mm

\subhead 4.3\endsubhead
Let $(\bb_{\g\n})$, $(\cb_{\g\n})$ be the bases of 
$\Gb\Ab_\g=\Hom_\AA(\Kc_\g,\AA)$ dual to $(IC_{\g\n})$, $(L_{\g\n})$.
Let $\otimes\,:\,\Gb\Ab_{\g'}\otimes\Gb\Ab_{\g''}\to\Gb\Ab_{\g'+\g''}$ 
and $\theta\,:\,\Gb\Ab_\g\to\Gb\Ab_\g$ 
be the maps dual to $\res_{\g'\g''}$ and $\DD.$
We consider the inductive system of $\AA$-modules 
$(\Gb\Ab_\g)$ such that $\bb_{\g\n}\mapsto\bb_{\g+\d,\n+\d}$.
Let $\Gb\Ab={\ds\lim_{\longrightarrow}}_\g\Gb\Ab_\g$ be the limit.
Let $\bb_\g,\cb_\g\in\Gb\Ab$ be the images of the elements
$\bb_{\g0}, \cb_{\g0}\in\Gb\Ab_\g$.

\proclaim{Theorem}
The $\AA$-module $\GRb$ is a subring of $\Kb(\Cc_q)$.
The linear map such that $\bb_\g\mapsto V(\g)$ 
is an algebra isomorphism $\Gb\Ab\simto\GRb$.
The map $\theta$ is a skew-linear antihomomorphism of $\Gb\Ab$ 
fixing the bases $\Bb=(\bb_\g)$, $\Cb=(\cb_\g)$.
For any $\g,\g'$ we have
$$\bb_\g\otimes\bb_{\g'}\in\bigoplus_{\g''}\NN[v^{-1},v]\cdot\bb_{\g''}.$$
\endproclaim

\noindent{\sl Proof :}
The maps $\DD$, $\res_{\g'\g''}$ are compatible 
with the projective system $(\Kc_\g)$. 
The limit, denoted $(\Kc,\res)$, is a co-algebra with
a skew-linear involution $\DD$.
By Lemma 4.2.1.$(a),(b)$, 
the projective system maps $IC_{\g\n}$ to $IC_{\g-\d,\n-\d}$,
for any $\n\in\Wedge^+(\g)$ such that $\n\succeq\d$.
In $\Kc$ we consider the elements 
$IC_\g=\bigl(IC_{\g+\d,\d}\bigr)$, with $\g\in X^+$,
and $L_{\g}=\bigl(L_{\g+\d,\d}\bigr)$, with $\g\in X$.
We have
$$\res(L_\g)=\sum_{\g=\g'+\g''}v^{\la\g',\g''\ra}L_{\g'}\otimes L_{\g''}.$$
Let $\Ab_X^\vee$ be the $\AA$-coalgebra with
the $\AA$-basis $(\ab_\g)$, $\g\in X$, and the coproduct
$\ab_\g\mapsto\sum_{\g=\g'+\g''}v^{\la\g',\g''\ra}\ab_{\g'}\otimes\ab_{\g''}.$
The $\AA$-linear map $\Ab_X^\vee\to(\Kc,\res)$ such that
$\ab_\g\mapsto L_\g$ is a surjective co-algebra homomorphism.
By Lemma 4.2.1.$(a)$ we have
$$L_\g=\sum_{\g'\in X^+}\gdm V(\g')_\g\cdot IC_{\g'},\quad\forall\g\in X.$$
The elements $IC_{\g'}$, $\g'\in X^+$, form a $\AA$-basis of $\Kc$.
Thus, the linear map 
$$\psi\,:\,\Gb\Ab\to\Ab_X,\quad
\bb_{\g'}\mapsto\sum_{\g\in X}\gdm V(\g')_\g\cdot e^\g,\quad
\forall\gamma'\in X^+,$$ 
is an injective ring homomorphism.
Consider the linear map $\phi\,:\,\Gb\Rb\to\Gb\Ab$ 
such that $V(\g)\mapsto\bb_\g$ for all $\g\in X^+$.
We get the commutative square of linear maps
$$\matrix
\Gb\Rb&{\buildrel\phi\over\longrightarrow}&\Gb\Ab\cr
\downarrow&&\downarrow{\ss\psi}\cr
\Kb(\Cc_q)&{\buildrel\gch\over\longrightarrow}&\Ab_X
\endmatrix$$
where $\psi,\gch$ are ring homomorphisms, see Proposition 2.4.$(a)$,
the vertical maps are injective, and $\phi$ is invertible.
Thus, $\Gb\Rb$ is a subring of $\Kb(\Cc_q)$ and $\phi$ is a ring homomorphism.
If $\g'+\g''=\g$ in $X^+$, then
$$(\DD\otimes\DD)\circ\res_{\g'\g''}\circ\DD=\res_{\g''\g'}.$$
Thus $\theta$ is an antihomomorphism.
\qed

\vskip3mm

If $\bb_\g\otimes\bb_{\g'}=v^{\la\g,\g'\ra}\bb_{\g+\g'}$, then
the $\Ub$-module $V(\g)\otimes V(\g')$ is simple and isomorphic to $V(\g+\g')$.
Conversely, if $V(\g)\otimes V(\g')$ is a simple $\Ub$-module 
it is isomorphic to $V(\g+\g')$. Then, the positivity in Theorem 4.3
implies that $\bb_\g\otimes\bb_{\g'}\in v^\ZZ\bb_{\g+\g'}$.
Then, by (2.3.1) we get
$\bb_\g\otimes\bb_{\g'}=v^{\la\g,\g'\ra}\bb_{\g+\g'}$.
The following conjecture generalizes to all
simply laced types the conjecture in \cite{BZ} (for type $A$).

\proclaim{Conjecture} The following statements are equivalent :

$\bb_\g\otimes\bb_{\g'}\in v^\ZZ\Bb,$
$\bb_\g\otimes\bb_{\g'}\in v^\ZZ\bb_{\g'}\otimes\bb_\g,$
and $\bb_\g\otimes\bb_{\g'}=v^{\la\g,\g'\ra}\bb_{\g+\g'}.$\hfill
\endproclaim


\head 5. The classical case\endhead
\subhead 5.1\endsubhead
Fix $\l,\l'\in P^+$.
Let $\Dc(N_{\l})^\heartsuit$,
$\Dc(N_{\l}\times N_{\l'})^\heartsuit$ 
be the full subcategories of 
$\Dc(N_\l)$,
$\Dc(N_{\l}\times N_{\l'})$ 
consisting of all complexes which are constructible 
with respect to the stratification in \S 3.2. 
Set $IC_{\l\a}=IC(N_{\l\a}^\heartsuit)$,
$\CC_{\l\a}=\CC_{N_{\l\a}^\heartsuit}[d_{\l\a}]$ if $\a\in\Wedge^+(\l)$, 
and set $L_{\l\a}=\pi_{\l\a!}\CC_{Q_{\l\a}}[d_{\l\a}]$ 
if $\a\in\Wedge(\l)$. 
Let $\Pc_\l$, $\Pc_{\l\l'}$ 
be the full subcategories of 
$\Dc(N_{\l})^\heartsuit$,
$\Dc(N_{\l}\times N_{\l'})^\heartsuit$
consisting of all complexes which are isomorphic to finite direct sums 
of complexes of the form $IC_{\l\a}$,
$IC_{\l\a}\boxtimes IC_{\l'\a'}$.

Assume that $\l=\l'+\l''$.
Setting $k=0$, $\g=\Id_W$ in \S 3.5 we get the commutative diagram
$$\matrix
Q_\l&{\buildrel\tilde\iota_\pm\over\hookleftarrow}&Q_\l^{\pm\t}
&{\buildrel\tilde\kappa_\pm\over\twoheadrightarrow}&
Q^\t_\l&\simeq&Q_{\l'}\times Q_{\l''}\cr
\downarrow&$\qed$&\downarrow&&\downarrow&&\downarrow\cr
N_\l&{\buildrel\iota_\pm\over\hookleftarrow}&N_\l^{\pm\t}
&{\buildrel\kappa_\pm\over\twoheadrightarrow}&N_\l^\t&
{\buildrel\phi\over\leftarrow}&N_{\l'}\times N_{\l''}.
\endmatrix$$
The restriction of the map $\tilde\kappa_\pm$ 
to $\tilde\kappa_\pm ^{-1}(Q_{\l'\a'}\times Q_{\l''\a''})$ is 
a vector bundle of rank
$$(d_{\l\a}-d_{\l'\a'}-d_{\l''\a''})/2,\leqno(1)$$
where $\a=\a'+\a''.$
Indeed, let $T_{\l\t}$ be the normal bundle to $Q_\l^\t$ in $Q_\l$,
and let $T_{\l\t}^\pm$ be the restriction to $Q_\l^\t$
of the relative tangent bundle to the map $\tilde\kappa_\pm$. 
The cocharacter $\t$ acts on $T_{\l\t}$ with non zero weights,
and $T_{\l\t}^\pm$ is the subbundle consisting 
of the positive (resp. negative) weights subspaces. 
Recall that $Q_\l$ has a $G_\l$-invariant holomorphic 
symplectic form, see \cite{N1, (3.3)}. 
Thus, the subvariety $Q_\l^\t$ is symplectic, and
the rank of $T_{\l\t}$ is twice the rank of $T_{\l\t}^\pm$.

Consider the functor
$$\res^\pm_{\l'\l''}=\kappa_{\pm!}\iota^*_\pm\,:\,
\Dc(N_{\l})^\heartsuit\to\Dc(N_\l^\t)^\heartsuit.$$
For any $\mu\in P^+$ we set
$$V(\l',\l'')_\mu=\Hom_{\gen}\bigl(V(\mu),V(\l')\otimes V(\l'')\bigr).$$

\proclaim{Lemma 1}
For any $\a\in\Wedge(\l)$ we have

(a) $\res^+_{\l'\l''}(L_{\l\a})=
\res^-_{\l'\l''}(L_{\l\a})=
\bigoplus_{\a=\a'+\a''}\phi_!(L_{\l'\a'}\boxtimes L_{\l''\a''}),$

(b) $\res_{\l'\l''}^\pm(IC_{\l\a})\simeq\Oplus_{\a',\a''} 
V(\l'-\a',\l''-\a'')_{\l-\a}\otimes
\phi_!(IC_{\l'\a'}\boxtimes IC_{\l''\a''}),$

(c) $\res^\pm_{\l'\l''}$ commutes to the Verdier duality.

(d) For any complex $P\in\roman{Ob}(\Pc_{\l})$ there is a complex
$P'\in\roman{Ob}(\Pc_{\l'\l''})$ such that
$\res^\pm_{\l'\l''}(P)\simeq\phi_!(P')$.
\endproclaim

\noindent{\sl Proof:}
Claim $(a)$ is proved as Lemma 4.1, using (5.1.1).
Using \cite{N2, Theorem 15.3.2} we get an isomorphism
$$L_{\l\a}\simeq
\Oplus_{\b\in Q^+}H_{top}(F_{\l-\b,\a-\b})\otimes IC_{\l\b}.\leqno (2)$$
Using Part $(a)$ and (5.1.2) we get
$$\Oplus_{\a\geq\b}V(\l-\b)_{\l-\a}\otimes\res^\pm_{\l'\l''}(IC_{\l\b})\simeq$$
$$\simeq\Oplus_{\a\geq\b}\Oplus_{\b',\b''}V(\l-\b)_{\l-\a}
\otimes V(\l'-\b',\l''-\b'')_{\l-\b}\otimes
\phi_!(IC_{\l'\b'}\boxtimes IC_{\l''\b''}),$$
where the sum is over all $\b',\b''\in Q^+$.
An induction on $\b$ gives
$$\res_{\l'\l''}^\pm(IC_{\l\b})\simeq
\Oplus_{\b',\b''} V(\l'-\b',\l''-\b'')_{\l-\b}\otimes
\phi_!(IC_{\l'\b'}\boxtimes IC_{\l''\b''}).$$
\qed

\vskip3mm

By Lemma 3.5.1.$(b)$, the functor $\phi_!$ is an equivalence from $\Pc_{\l'\l''}$
to a full subcategory of $\Dc(N_\l^\tau)^\heartsuit$.
Composing $\res_{\l'\l''}^\pm$ with a quasi-inverse to $\phi_!$ we get a functor
$\res_{\l'\l''}\,:\,\Pc_\l\to\Pc_{\l'\l''}$.
Let $\Vc ec$ be the category of finite dimensional
complex vector spaces, 
and let $\Pc_\l^\circ$ be the category dual to $\Pc_\l$.
We consider the following functors
$$\matrix
\odot\,:\hfill&\Pc_{\l'}^\circ\times\Pc_{\l''}^\circ\to
\Pc_\l^\circ,\hfill&(P',P'')\mapsto\Oplus_\a
\Hom_{\Pc_{\l'\l''}}\bigl(\res_{\l'\l''}(IC_{\l\a}), 
P'\boxtimes P''\bigr)\otimes IC_{\l\a},\hfill\cr\cr
\Phi_\l\,:\hfill&\Pc_\l^\circ\to\Vc ec,\hfill&
P\mapsto\Hom_{\Pc_\l}\bigl(P,\Oplus_\a L_{\l\a}\bigr),\hfill\cr\cr
p_{\l\b}\,:\hfill&\Pc_\l^\circ\to\Pc_{\l-\b}^\circ,\hfill&
P\mapsto\Oplus_{\a\geq\b}\Hom_{\Pc_\l}
(IC_{\l\a},P)\otimes IC_{\l-\b,\a-\b},\hfill
\endmatrix$$
where $\b\in\Wedge^+(\l)$.
Note that \cite{N2, Theorem 3.3.2} and \S 3.2 give
$$\a-\b\in\Wedge^+(\l-\b)\iff\a\in\Wedge^+(\l),\quad\a\geq\b,$$
and similarly with $\Wedge(\l)$.
By (5.1.2) we have 
$$p_{\l\b}(L_{\l\a})\simeq\left\{\matrix
L_{\l-\b,\a-\b}\quad&\text{if}\quad\a\geq\b\hfill\cr
0\quad&\text{else}.\hfill
\endmatrix\right.$$

We define a new category $\Pc^\circ$ as follows.
Objects of $\Pc^\circ$ are collections $P=(P_\l,\g_{\l\b})$, 
where $\l\in P^+$, $\b\in\Wedge^+(\l)\setminus\{0\}$, 
$P_\l\in\roman{Ob}(\Pc_\l)$ and
$$\g_{\l\b}\in\Isom_{\Pc_{\l-\b}}
\bigl(P_{\l-\b},p_{\l\b}(P_\l)\bigr)$$ 
are isomorphisms satisfying the obvious chain condition. 
Morphisms $P'\to P''$ are 
collections $(\phi_\l)\in\Prod_\l\Hom_{\Pc_\l}(P''_\l,P'_\l)$
such that
$$\g'_{\l\b}\circ\phi_{\l-\b}=p_{\l\b}(\phi_\l)\circ\g''_{\l\b}
\in\Hom_{\Pc_{\l-\b}}\bigl(P''_{\l-\b},p_{\l\b}(P'_\l)\bigr).$$

\proclaim{Lemma 2} 
Fix $\b\in\Wedge^+(\l)$, 
$\b'\in\Wedge^+(\l')$, $\b''\in\Wedge^+(\l'')$ such that $\b=\b'+\b''$.
For any $P,P',P''\in\roman{Ob}(\Pc^\circ)$ we have natural embeddings
$$\Phi_{\l-\b}(P_{\l-\b})\subset\Phi_\l(P_\l),\quad
\quad P'_{\l'-\b'}\odot P''_{\l''-\b''}\subset
p_{\l\b}(P'_{\l'}\odot P''_{\l''}).$$
Moreover we have
$\Sum_{\b',\b''}P'_{\l'-\b'}\odot P''_{\l''-\b''}=
p_{\l\b}(P'_{\l'}\odot P''_{\l''}).$
\endproclaim

\noindent{\sl Proof:}
Fix an isomorphism as in (5.1.2) for each $\a\in\Wedge(\l)$.
For any such $\a$ we get a morphism of functors
$$\Oplus_{\a'}\Hom(-,IC_{\l\a'})\otimes\Hom(IC_{\l-\b,\a'-\b},L_{\l-\b,\a-\b})\to 
\Hom(-,L_{\l\a}).$$
By definition of $\Phi_\l$, $p_{\l\b}$ this morphism gives a morphism of functors
$\Phi_{\l-\b}\circ p_{\l\b}\to\Phi_\l$.
The morphism
$\Phi_{\l-\b}(P_{\l-\b})\to\Phi_\l(P_\l)$ is the composition
of the isomorphism $\Phi_{\l-\b}(\g_{\l\b})$ and the morphism 
$\Phi_{\l-\b}\circ p_{\l\b}\to\Phi_\l$ above.
Using (5.1.2) we get
$$\Phi_\l(IC_{\l\a})\simeq V(\l-\a),\quad
\Phi_{\l-\b}\circ p_{\l\b}(IC_{\l\a})\simeq\left\{
\matrix 
V(\l-\a)\quad&\roman{if}\quad\a\geq\b,\hfill\cr
0\quad&\roman{else}.\hfill\cr
\endmatrix\right.$$
This proves Claim one.
For any $\a'\in\Wedge^+(\l')$, $\a''\in\Wedge^+(\l'')$
Lemma 5.1.1 gives an isomorphism of complexes 
$$IC_{\l'\a'}\odot IC_{\l''\a''}\simeq
\Oplus_\a V(\l'-\a',\l''-\a'')_{\l-\a}\otimes IC_{\l\a},\leqno(3)$$
where the sum is over all $\a\in\Wedge^+(\l)$ such that
$V(\l'-\a',\l''-\a'')_{\l-\a}\neq\{0\}$. 
Fix such a family of isomorphisms.
It gives a morphism of functors
$p_{\l'\b'}(-)\odot p_{\l''\b''}(-)\to p_{\l\b}(-\odot -)$.
The morphism
$P'_{\l'-\b'}\odot P''_{\l''-\b''}\to p_{\l\b}(P'_{\l'}\odot P''_{\l''})$
is the composition of the isomorphism
$\g'_{\l'-\b'}\odot\g''_{\l''-\b''}$
and the morphism of functors  
$p_{\l'\b'}(-)\odot p_{\l''\b''}(-)\to p_{\l\b}(-\odot -)$
above. Then, Claim two and three are consequences of the following identities.
If $V(\l'-\a',\l''-\a'')_{\l-\a}\neq\{0\}$, then
$\a\geq\a'+\a''$, and thus
$$\matrix
\a\geq\b&\Leftarrow&\a'\geq\b',\a''\geq\b'',\hfill\cr\cr
\a\geq\b\hfill&\Rightarrow&\exists\,\b',\b''\quad\roman{s.t.}\quad
\a'\geq\b',\a''\geq\b'',\b=\b'+\b''.\hfill
\endmatrix$$
We are done.
\qed

\vskip3mm

By Lemma 5.1.2 the category $\Pc^\circ$ is endowed with the functors
$\Phi\,:\,\Pc^\circ\to\Vc ec$, 
$\odot\,:\,\Pc^\circ\times\Pc^\circ\to\Pc^\circ$
such that 
$$\Phi(P)={\lim_\lra}_\l\Phi_\l(P_\l),\quad
(P'\odot P'')_\l=\Sum_{\l=\l'+\l''}P'_{\l'}\odot P''_{\l''}.$$
Then, (5.1.3) gives  the following.

\proclaim{Lemma 3} 
$(\Pc^\circ,\odot)$ is a tensor category, and $\Phi$ is a tensor functor.
\endproclaim

\noindent
Let $\Ab$ be the Grothendieck group of $\Pc^\circ$.
The functor $\odot$ gives a product $\Ab\otimes\Ab\to\Ab.$
Let $\bb_\l$, $\cb_\l$ be the classes in $\Ab$ of 
the objects of $\Pc^\circ$ associated to the families
$(IC_{\l+\b,\b})$, $(L_{\l+\b,\b})$.
Then $(\bb_\l)$, $(\cb_\l)$ are bases of $\Ab$.
Let $(\Rc(\gen),\otimes)$ be
the tensor category of finite dimensional $\gen$-modules.
We have proved the following theorem.

\proclaim{Theorem}
The tensor categories $(\Pc^\circ,\odot)$,
$(\Rc(\gen),\otimes)$ are equivalent.
The group homomorphism such that $\bb_\l\mapsto V(\l)$ is a ring isomorphism
$\Ab\simto\Rb(\gen)$. 
Moreover, we have
$\bb_\l=\sum_{\mu}\dim V(\l)_{\mu}\cdot\cb_{\mu}$.
\endproclaim

\subhead 5.2\endsubhead
In this subsection we consider the non simply laced case.
Our construction is based on \cite{L1, \S 11}.
Assume that $\underline\gen$ is a non simply laced,
simple, complex Lie algebra.
Fix a simply laced simple Lie algebra $\gen$ and a diagram
automorphism $a$ of $\gen$ such that
the Dynkin graph of $\underline\gen$ is deduced
from the Dynkin graph of $\gen$ 
as in \cite{L1, \S 14}.
Let $n$ be the order of the automorphism $a$ 
($n=2$ for types $B_k, C_k, F_4$, and $n=3$ for type $G_2$).
The automorphism $a$ is identified with a permutation of the set 
$I\times H$, see \S 3.1, such that
$$a(h')=a(h)',\quad a(h'')=a(h)'',\quad 
a(\overline h)=\overline{a(h)}.$$
Let $\langle a\rangle$ be the cyclic group of automorphisms
of $(I,H)$ generated by $a$. 
Let $\Iu$ be the set of $\la a\ra$-orbits in $I,$ 
and let $\Pu^+=(P^+)^a$, $\Qu^+=(Q^+)^a$ 
be the corresponding sub-semigroups of $P^+$, $Q^+$.
The simple root $\a_\iu$ and the fundamental weight $\o_\iu$ of 
$\underline\gen$
are identified with the sums $\sum_{i\in\iu}\a_i\in\Qu^+$, 
$\sum_{i\in\iu}\o_\iu\in\Pu^+$. 
For any $\l\in P^+$, $\a\in Q^+$, 
the diagram automorphism induces natural isomorphisms 
$Q_{\l\a}\simto Q_{a(\l), a(\a)}$, $N_{\l\a}\simto N_{a(\l), a(\a)}$.
Let denote them by $a$ again.

To avoid confusions, finite dimensional representations of
$\gen$, $\underline\gen$
are denoted by $V(\l)$, $\Vu(\l)$ respectively.
The subsets of $Q^+$, $\Qu^+$ defined in \S 2.1 are denoted
by $\Wedge(\l)$, $\Wedge^+(\l)$ and
$\underline\Wedge(\l)$, $\underline\Wedge^+(\l)$ respectively.

Fix $\l,\l'\in\Pu^+$ and $\a\in Q^+$. 
Following \cite{L1, \S 11} we consider new categories $\aa\Pc_\l$,
$\aa\Pc_{\l\l'}$.
An object of $\aa\Pc_\l$ is a pair $(P,\theta)$, 
where $P\in\roman{Ob}(\Pc_\l)$ and 
$\theta\,:\,a^*P\simto P$
is an isomorphism such that the composition
$$a^{*n}P\lra\cdots\lra a^{*2}P{\buildrel a^*\theta\over\lra}a^*P
{\buildrel\theta\over\lra}P$$
is the identity.
A morphism $(P,\theta)\to(P',\theta')$ is a morphism $f\,:\,P\to P'$
such that $f\theta=\theta'(a^*f)$.
The category $\aa\Pc_{\l\l'}$ is constructed in the same way.
Both categories are Abelian.
For any functor $F\,:\,\Pc_\l\to\Pc_{\l'}$
and for any isomorphism of functor $a^*F\simto Fa^*$
there is the functor $\aa F\,:\,\aa\Pc_\l\to\aa\Pc_{\l'}$
such that $\aa F(P,\theta)=(F(P),\theta^F)$ where $\theta^F$ 
is the composition of the chain of maps
$$a^*F(P)\longrightarrow 
F(a^*P){\buildrel F(\theta)\over\longrightarrow}F(P).$$

The functor $a^*$ on $\Pc_\l$ has the order $n$, where $n=2$ or 3.
Let $\aa\Ic_\l$ be the full subcategory of $\aa\Pc_\l$ whose objects
are the pairs $(P,\theta)$ such that 
$P\simeq P'\oplus a^*P'\oplus\cdots\oplus(a^*)^{n-1}P'$
for some $P'\in\Pc_\l$, and $\theta$ is an isomorphism
carrying the direct summand $(a^*)^jP'\subset a^*P$ onto the direct summand
$(a^*)^jP'\subset P$. 
The objects of $\aa\Ic_\l$ are said to be traceless.

The automorphism $a$ preserves the stratification of $N_\l$.
Since $IC_{\l\a}$ is canonically attached to $N_{\l\a}^\heartsuit$, 
there is a canonical isomorphism
$a^*IC_{\l,a(\a)}\simto IC_{\l\a}.$ 
If $\a\in\Qu^+$ the corresponding object in $\aa\Pc_\l$ 
is denoted by $\aa IC_{\l\a}.$ 
Let $\mu_n\subset\CC^\times$ be the set of $n$-th roots of unity.
For any $\zeta\in\mu_n$ and any $Q=(P,\theta)\in\roman{Ob}(\aa\Pc_\l)$
we put $Q(\zeta)=(P,\zeta\theta)$.
If $\a\notin\Qu^+$ and $\zeta_1,...,\zeta_n\in\mu_n$, let
$\aa IC_{\l\a}(\zeta_1,...,\zeta_n)$ be the object of $\aa\Pc_\l$ 
associated to the perverse sheaf 
$$P=IC_{\l\a}\oplus IC_{\l,a(\a)}\oplus\cdots IC_{\l,a^{n-1}(\a)}$$
and the isomorphism $a^*P\simto P$ which maps the summand
$a^*IC_{\l,a^i(\a)}$ onto the summand $IC_{\l,a^{i-1}(\a)}$
by $\zeta_{i+1}$ times the canonical isomorphism.
A simple object in $\aa\Pc_\l$ is isomorphic either to
$\aa IC_{\l\a}(\zeta)$ for some $\a\in\Qu^+$ and $\zeta\in\mu_n$, 
or to $\aa IC_{\l\a}(\zeta_1,...,\zeta_n)$ for some
$\a\in Q^+\setminus\Qu^+$ and $\zeta_1,...,\zeta_n\in\mu_n$.
Let $\11\Pc_\l$ be the full subcategory of $\aa\Pc_\l$ whose objects
are isomorphic to finite direct sums of the objects $\aa IC_{\l\a}$.

The image by the functor $\pi_{\l\a!}$
of the obvious isomorphism $a^*\CC_{Q_{\l,a(\a)}}\simto\CC_{Q_{\l\a}}$
is an isomorphism $a^*L_{\l,a(\a)}\simto L_{\l\a}.$ 
If $\a\in\Qu^+$ the corresponding object in $\aa\Pc_\l$ 
is denoted by $\aa L_{\l\a}$.
Assume that $\b\in Q^+$ is such that 
$\a\geq\b$ and $N_{\l\b}^\heartsuit\neq\emptyset$.
Fix an element $x_\b\in N_{\l\b}^\heartsuit$.
One proves as in \cite{N2, Theorem 3.3.2} that there are
$\la a\ra$-invariant open sets 
$$U_\a\subset\la a\ra(N_{\l\a}),\quad
U^\heartsuit_\b\subset\la a\ra(N^\heartsuit_{\l\b}),\quad
U_{\a-\b}\subset\la a\ra(N_{\l-\b,\a-\b})$$
containing $x_\b$, $x_\b$, 0 respectively, 
and a commutative square 
$$\matrix
U_\a&\simto&U^\heartsuit_\b\times U_{\a-\b}\cr
{\ss\pi}\uparrow&&\uparrow{\ss\Id\times\pi}\cr
\pi^{-1}(U_\a)&\simto&
U^\heartsuit_\b\times\pi^{-1}(U_{\a-\b}),
\endmatrix$$
where $\pi$ denotes either $\pi_{\l\a}$ or $\pi_{\l-\b,\a-\b}$.
The horizontal maps are analytic $\la a\ra$-equivariant isomorphisms
carrying the element $x_\b\in U_\a$ to 
$(x_\b,0)\in U_\b^\heartsuit\times U_{\a-\b}$.
By (5.1.2) we have
$$L_{\l\a}\simeq\Oplus_{\b\in Q^+}
H_{top}(F_{\l-\b,\a-\b})\otimes IC_{\l\b},$$
and the isomorphism $a^*L_{\l,a(\a)}\simto L_{\l\a}$ 
maps the direct summand
$$H_{top}\bigl(a(F_{\l-\b,\a-\b})\bigr)\otimes a^* IC_{\l,a(\b)}
\quad\roman{onto}\quad 
H_{top}(F_{\l-\b,\a-\b})\otimes IC_{\l\b}$$
in the obvious way. 
By \cite{X, Theorem 3.2.1}, if $\a,\b\in\Qu^+$
the number of irreducible components
of $F_{\l-\b,\a-\b}$
which are mapped to themselves by $a$
is the multiplicity $\dim\Vu(\l-\b)_{\l-\a}$.
Thus $\aa L_{\l\a}=\11 L_{\l\a}\oplus I_{\l\a}$ where
$$\11 L_{\l\a}\simeq
\Oplus_{\b\in\Qu^+}\Vu(\l-\b)_{\l-\a}\otimes\aa IC_{\l\b}
\in\roman{Ob}(\11\Pc_\l),
\quad I_{\l\a}\in\roman{Ob}(\aa\Ic_\l).\leqno(1)$$

Assume that $\l=\l'+\l''$ in $\Pu^+$.
The maps $\iota_\pm$, $\kappa_\pm$, $\phi$
commute to the automorphism $a$ of $N_\l$.
Thus, there is a natural isomorphism
$a^*\res_{\l'\l''}\simto\res_{\l'\l''}a^*.$
We get the functor
$\aa\res_{\l'\l''}\,:\,\aa\Pc_\l\to\aa\Pc_{\l'\l''}$.
Lemma 5.1.1 implies the following.

\proclaim{Lemma}
For any $\a\in\Qu^+$  there are traceless objects $I,I'$ such that

(a) $\aa\res_{\l'\l''}(\aa L_{\l\a})=I\oplus
\bigoplus_{\a=\a'+\a''}\aa L_{\l'\a'}\boxtimes\aa L_{\l''\a''},$

(b) $\aa\res_{\l'\l''}(\aa IC_{\l\a})=I'\oplus\Oplus_{\a',\a''} 
\Vu(\l'-\a',\l''-\a'')_{\l-\a}\otimes 
(\aa IC_{\l'\a'}\boxtimes\aa IC_{\l''\a''}).$
\endproclaim

For any $\b\in\Wedge^+(\l)\cap\Qu^+$, 
there is an obvious isomorphism of functors
$a^*p_{\l\b}\simto p_{\l\b}a^*$.
The corresponding functor 
$\aa p_{\l\b}\,:\,\aa\Pc_\l\to\aa\Pc_{\l-\b}$ 
is exact and satisfies
$$\aa p_{\l\b}(\aa IC_{\l\a})=\left\{\matrix
\aa IC_{\l-\b,\a-\b}\quad&\text{if}\quad\a\geq\b\hfill\cr
0\quad&\text{else}.\hfill
\endmatrix\right.$$
We have also, see (5.2.1),
$$\aa p_{\l\b}(\aa L_{\l\a})\simeq\left\{\matrix
\aa L_{\l-\b,\a-\b}\quad&\text{if}\quad\a\geq\b\hfill\cr
0\quad&\text{else}.\hfill
\endmatrix\right.$$
Let $\Kb(\aa\Pc_\l)$ be the Grothendieck group of $\aa\Pc_\l$.
The class of an object $P$ is still denoted by $P$.
Let $\kb\subset\CC$ be subring generated by $\mu_n$.
Let $\Kc'_\l$ be the quotient of $\Kb(\aa\Pc_\l)\otimes\kb$
by the relations :

\itemitem{--} $P(\zeta)=P\otimes\zeta$ for any $\zeta\in\mu_n$,

\itemitem{--} the class of a traceless object is zero.

\noindent
Let $\Kc_\l\subset\Kc'_\l$ be the subgroup spanned by the classes
of objects in $\11\Pc_\l$, and let
$\Ab_\l=\Kc_\l^*$ be the dual group.
Using the maps $\aa p_{\l\b}$ we construct, as in \S 4.3,
an inductive system of groups $(\Ab_\l)$.
The limit, denoted by $\Ab$, is endowed
with a product $\Ab\otimes\Ab\to\Ab,$
and two distinguished bases
$(\bb_\l)$, $(\cb_\l)$ 
associated to the families
$(\aa IC_{\l+\b,\b})$, $(\aa L_{\l+\b,\b})$.

\proclaim{Theorem}
The group homomorphism such that 
$\bb_\l\mapsto\Vu(\l)$ is a ring isomorphism 
$\Ab\simto\Rb(\underline\gen)$. 
Moreover, we have
$\bb_\l=\sum_{\mu}\dim\Vu(\l)_{\mu}\cdot\cb_{\mu}$.
\endproclaim

\subhead 5.3\endsubhead
In this subsection we explain how a similar construction gives a 
natural restriction map $\Gb\Rb\to\Rb(\gen)\otimes\AA$.
Consider the diagram
$$\matrix
Q_\l&{\buildrel\tilde\iota_\pm\over\hookleftarrow}&Q_\l^{\pm\g}
&{\buildrel\tilde\kappa_\pm\over\twoheadrightarrow}&Q_\g\cr
\downarrow&$\qed$&\downarrow&&\downarrow\cr
N_\l&{\buildrel\iota_\pm\over\hookleftarrow}&N_\l^{\pm\g}
&{\buildrel\kappa_\pm\over\twoheadrightarrow}&N_\g.
\endmatrix$$
Set $\eps_\g=\eps_{\g\g}$. 
Let $\kappa^\pm_\n$ be the relative dimension of $\tilde\kappa_\pm$ 
above the component $Q_{\g\n}$.
The same computations as in Lemma 3.5.2 or in (5.1.1) give
$$\matrix
\kappa^-_\n=d_{\l\a}/2-d_{\g\n},\hfill&\kappa^+_\n=d_{\l\a}/2,\hfill\cr
\eps_\g-\eps_{\g-\d}=d_{\g\d},\hfill&d_{\g\n}-d_{\g-\d,\n-\d}=d_{\g\d},\hfill
\endmatrix$$
for any $\d\in\Wedge^+(\g)$, $\n\succeq\d.$
Consider the functor
$$\res^\pm_\g=\kappa_{\pm!}\iota^*
\,:\,\Dc(N_\l)^\heartsuit\to \Dc(N_\g)^\heartsuit.$$
By base change we get, for any $\n\in\Wedge^+(\g)$,
$$\res^\pm_\g L_{\l\a}=\pi_{\g!}\tilde\kappa_{\pm!}\tilde\iota_\pm^*
\CC_{Q_{\l\a}}[d_{\l\a}]
=\Oplus_\n\pi_{\g!}\CC_{Q_{\g\n}}[d_{\l\a}-2\kappa^\pm_\n]=
\Oplus_\n L_{\g\n}[\mp d_{\g\n}].$$

\proclaim{Lemma} 
(a) For any complex $P\in\roman{Ob}(\Pc_\l)$ the complex
$\res^\pm_\g(P)$ belongs to $\roman{Ob}(\Pc_\g)$. 

(b) We have $\DD\circ\res^+_\g=\res^-_\g\circ\DD.$
\endproclaim

\noindent
The corresponding group homomorphism 
$v^{\eps_\g}\res^+_\g\,:\,\Kb(\Pc_\l)\to\Kc_\g$ 
is compatible with the projective systems in \S 5.1, \S 4.2. 
Let $\res\,:\,\Gb\Ab\to\Ab\otimes\AA$ 
be the inductive limit of the system of maps dual to
$v^{\eps_\g}\res^+_\g.$ 

\proclaim{Proposition}
The element $\res(\bb_\g)$ belongs to $\Oplus_\l\NN[v^{-1},v]\cdot\bb_\l$ 
for all $\g\in X^+$. If $\g\in G_\l^{\vee,\ad}$ then
$\res(\cb_\g)=v^{\eps_\g}\cb_\l$. 
\endproclaim

\vfill
{\eightpoint{
$$\matrix\format\l&\l&\l&\l\\
\phantom{.} & {\text{Michela Varagnolo}}\phantom{xxxxxxxxxxxxx} &
{\text{Eric Vasserot}}\\
\phantom{.}&{\text{D\'epartement de Math\'ematiques}}\phantom{xxxxxxxxxxxxx} &
{\text{D\'epartement de Math\'ematiques}}\\
\phantom{.}&{\text{Universit\'e de Cergy-Pontoise}}\phantom{xxxxxxxxxxxxx} &
{\text{Universit\'e de Cergy-Pontoise}}\\
\phantom{.}&{\text{2 Av. A. Chauvin}}\phantom{xxxxxxxxxxxxx} & 
{\text{2 Av. A. Chauvin}}\\
\phantom{.}&{\text{95302 Cergy-Pontoise Cedex}}\phantom{xxxxxxxxxxxxx} & 
{\text{95302 Cergy-Pontoise Cedex}}\\
\phantom{.}&{\text{France}}\phantom{xxxxxxxxxxxxx} & 
{\roman{France}}\\
&{\text{email: michela.varagnolo\@u-cergy.fr}}\phantom{xxxxxxxxxxxxx} &
{\text{email: eric.vasserot\@u-cergy.fr}}
\endmatrix$$
}}

\vskip1cm
\Refs
\widestnumber\key{ABC}

\ref\key{BZ}\by{\it Berenstein, A., Zelevinsky, A.}
\paper{\rm String bases for quantum groups of type $A_r$}
\jour Advances in Soviet Math.\vol 16\yr 1993\pages 51-89\endref

\ref\key{CG}\by{\it Chriss, N., Ginzburg, V.}
\book{\rm Representation theory and complex geometry}
\publ Birkh\"auser\publaddr Boston-Basel-Berlin\yr 1997\endref

\ref\key{D}\by{\it Damiani, I.}
\paper{\rm La R-matrice pour les alg\`ebres quantiques de type affine non tordu}
\jour Ann. Sci. \'Ecole Norm. Sup. (4)\vol 31\yr 1998
\pages 493-523\endref

\ref\key{FM}\by{\it Frenkel, E., Mukhin, E.}
\paper{\rm Combinatorics of $q$-characters of finite-dimensional
representations of quantum affine algebras}\jour math.QA/9911112\endref

\ref\key{FR}\by{\it Frenkel, E., Reshetikhin, N.}
\paper{\rm The $q$-characters of representations of quantum affine algebras
and deformations of $W$-algebras}
\jour Contemporary Math.\vol 248\yr 2000\pages 163-205
\endref

\ref\key{GV}\by{\it Ginzburg, V., Vasserot, E.}
\paper{\rm Langlands reciprocity fo affine quantum groups of type $A_n$}
\jour Internat. Math. Res. Notices\vol 3\yr 1993\pages 67-85\endref

\ref\key{Ka}\by{\it Kashiwara, M.}
\paper{\rm On level zero representations of quantized affine algebras}
\jour math.QA/0010293\endref

\ref\key{Kn}\by{\it Knight, H.}
\paper{\rm Spectra of tensor products of finite-dimensional representations
of Yangians}
\jour J. Algebra\vol 174\yr 1994\pages 187-196
\endref

\ref\key{L1}\by{\it Lusztig, G.}
\book{\rm Introduction to quantum groups}
\publ Birkh\"auser\publaddr Boston-Basel-Berlin \yr 1994\endref

\ref\key{L2}\by{\it Lusztig, G.}
\paper{\rm On quiver varieties}
\jour Adv. in Math.\vol 136\yr 1998\pages 141-182\endref

\ref\key{LNT}\by{\it Leclerc, B., Nazarov, M., Thibon, J.-Y.}
\paper{\rm Induced representations of affine Hecke algebras and canonical bases
of quantum groups}
\jour math.QA/0011074\endref

\ref\key{M}\by{\it Malkin, A.}
\paper{\rm Tensor product varieties and crystals. ADE case}
\jour math.AG/0103025\endref

\ref\key{N1}\by{\it Nakajima, H.}
\paper{\rm Quiver varieties and Kac-Moody algebras}
\jour Duke Math. J.\vol 91\yr 1998\pages 515-560\endref

\ref\key{N2}\by{\it Nakajima, H.}
\paper{\rm Quiver varieties and finite dimensional 
representations of quantum affine algebras}
\jour Jour. A. M. S.\vol 14\yr 2001\pages 145-238\endref

\ref\key{N3}\by{\it Nakajima, H.}
\paper{\rm $t$-analogue of the $q$-character
of finite dimensional representations of quantum affine algebras}
\jour math.QA/0009231\endref

\ref\key{Vr}\by{\it Varagnolo, M.}
\paper{\rm Quiver varieties and Yangians}
\jour Letters in Math. Phys.\vol 53\yr 2000\pages 273-283\endref

\ref\key{Va}\by{\it Vasserot, E.}
\paper{\rm Affine quantum groups and equivariant K-theory}
\jour Transformation groups\vol 3\yr 1998\pages 269-299\endref

\ref\key{VV}\by{\it Varagnolo, M., Vasserot, E.}
\paper{\rm Standard modules of quantum affine algebras}
\jour math.QA/0006084\endref

\ref\key{X}\by{\it Xu, F.}
\paper{\rm A note on quivers with symmetries}
\jour math.QA/9707003\endref
\endRefs
\enddocument